\documentclass[11pt,leqno]{article}
\normalfont
\renewcommand{\baselinestretch}{1.24}
\usepackage{amsfonts}
\usepackage{mathrsfs}
\usepackage{esint}
\usepackage{graphicx}

\newfont{\eulercursive}{eurm10 at 11pt}

\newcommand{\myr}{\mbox{\eulercursive r}}
\newcommand{\mys}{\mbox{\eulercursive s}}
\newcommand{\myeulerR}{\mbox{\eulercursive R}}

\newfont{\smallereulercursive}{eurm10 at 10pt}
\newcommand{\mynewd}{\mbox{\smallereulercursive d}}
\newcommand{\mynewr}{\mbox{\smallereulercursive r}}

\newfont{\muchsmallereulercursive}{eurm10 at 8pt}
\newcommand{\mynewerd}{\mbox{\muchsmallereulercursive d}}
\newcommand{\mynewerr}{\mbox{\muchsmallereulercursive r}}
\newcommand{\mynewers}{\mbox{\muchsmallereulercursive s}}

\newcommand{\QED}{\raisebox{0.5mm}{\fbox{\rule{0mm}{1.5mm}\ }}}

\newcommand{\canonicalU}{\mbox{\large $\mathscr{U}\!$}}

\newcommand{\myA}{\mbox{\sffamily A}}

\newcommand{\mysmallerA}{\mbox{\scriptsize \sffamily A}}

\newcommand{\myP}{\mbox{\sffamily P}}

\newcommand{\myS}{\mbox{\sffamily S}}

\newcommand{\mysmallerS}{\mbox{\scriptsize \sffamily S}}

\newcommand{\myR}{\mbox{\sffamily R}}

\newcommand{\mysmallerR}{\mbox{\scriptsize \sffamily R}}

\begin{document}

\newpage
\setcounter{page}{1} 
\renewcommand{\baselinestretch}{1}

\vspace*{-0.7in}
\hfill {\footnotesize February 19, 2019 (a minor variation of a July 25, 2018 edition)}

\begin{center}
{\Large \bf Unit shapes and a wealth of calculus problems}

\vspace*{0.05in}
\renewcommand{\thefootnote}{1}
Robert G.\ Donnelly\footnote{Department of Mathematics and Statistics, Murray State
University, Murray, KY 42071\\ 
\hspace*{0.25in}Email: {\tt rob.donnelly@murraystate.edu}} 
\renewcommand{\thefootnote}{2} 
\hspace*{-0.07in}and Alexander F.\ Thome\footnote{Department of Mathematics and Statistics, Murray State
University, Murray, KY 42071\\ 
\hspace*{0.25in}Email: {\tt athome1@murraystate.edu}}

\end{center}

\begin{abstract}
For a given family of similar shapes, what we call a ``unit shape'' strongly analogizes the role of the unit circle within the family of all circles. 
Within many such families of similar shapes, we present what we believe is naturally and intrinsically unital about their unit shapes. 
We present a number of calculus problems related to extremal questions about collections of unit shapes, and we recapitulate some isoperimetric problems in terms of unit shapes. 
We close by presenting some problems (some of which are open) and by proffering perhaps a new perspective on the $\pi$ vs.\ $\tau$ debate. 

\begin{center}


\end{center}
\end{abstract}

{\bf \S 1\ \ Introduction.} 
A striking and frequently pondered fact commonly encountered in first-semester Calculus is that 
\[\frac{d}{dr}\left(\pi  r^{2}\right) = \pi \cdot 2r,\]
an equality which can be interpreted as saying that, for circles, the derivative of area is perimeter (i.e.\ circumference) when each of these quantities is understood with respect to circle radius $r$. 
This familiar identity is often interpreted/explained using differentials, as one can see, for example, in the internet discussions \cite{MSE} and \cite{AP}. 
Here is how this line of reasoning generally goes: Express circle area $\myA$ and perimeter $\myP$ as functions of radius: $\myA(r)=\pi\, r^{2}$ and $\myP(r) = \pi \cdot 2r$. 
Then a small change $dr$ in the radius $r$ results in an exact change in area 
\[\Delta\myA = \pi(r+dr)^{2}-\pi r^{2} = \pi \cdot 2\left[\frac{r+(r+dr)}{2}\right]dr = \pi \cdot 2rdr + \pi(dr)^{2},\] 
which is the area of an annulus of thickness $dr$ along the circle's boundary. 
If $dr$ is small, then this thin piece of area is about the same as the area of a rectangle of length $\pi \cdot 2r$ and height $dr$, i.e.\ the contribution of the term $\pi(dr)^{2}$ is relatively negligible. 
So $\Delta\myA \approx d\myA = \pi \cdot 2rdr = \myP(r)dr$ as $dr \rightarrow 0$. 
One could view this geometric reasoning as an interpretation/explanation of the differential relationship   
\begin{equation}
\myA'(r) = \frac{d\myA}{dr} = \myP(r)  
\end{equation}
for circles. 

Perhaps one reason that fact (1) catches our eye is that it does not immediately square with other similar scenarios. 
For example, if we express circle area and perimeter in terms of diameter $\mynewd$, then $\myA(\mynewd) = \frac{\pi}{4}\mynewd^{2}$, $\myP(\mynewd) = \pi\mynewd$, but $\frac{d\mysmallerA}{d\mynewerd} = \myA'(\mynewd) = \frac{\pi}{2}\mynewd \ne \pi\mynewd = \myP(\mynewd)$. 
Or, if we express square area $\myA$ and perimeter $\myP$ in terms of side length $s$, then $\myA(s) = s^{2}$, $\myP(s) = 4s$, but $\frac{d\mysmallerA}{ds} = \myA'(s) = 2s \ne 4s = \myP(s)$. 

Among the many papers that take more general viewpoints on fact (1) --- see for example \cite{BYU}, \cite{NCTM} ---  
our favorite is from a {\em College Mathematics Journal} note from 1997 by Jingcheng Tong \cite{Tong}. 
Let us briefly summarize the perspective of \cite{Tong} for planar figures using notions that we more fully develop in \S 2 and \S 3 below. 
Let $\mathcal{C}$ be a shape in the plane with perimeter $\myP(\mathcal{C})$, semiperimeter $\myS(\mathcal{C}) = \frac{1}{2}\myP(\mathcal{C})$, and area $\myA(\mathcal{C})$. 
Suppose $\lambda$ is a positive real number. 
Apply the planar scaling function that stretches or compresses the plane by a factor of $\lambda$. 
Let $\mathcal{C}_{\lambda}$ be the image of shape $\mathcal{C}$ under this transformation. 
In \S 2 we will more carefully consider the following intuitively plausible facts: $\myP(\mathcal{C}_{\lambda}) = \lambda\myP(\mathcal{C})$, $\myS(\mathcal{C}_{\lambda}) = \lambda\myS(\mathcal{C})$, and $\myA(\mathcal{C}_{\lambda}) = \lambda^{2}\myA(\mathcal{C})$. 
Next, set $\upsilon := \myS(\mathcal{C})/\myA(\mathcal{C})$. 
If we let $\mathcal{D}_{\lambda} := \mathcal{C}_{\upsilon\lambda}$, then $\lambda = \myA(\mathcal{D}_{\lambda})/\myS(\mathcal{D}_{\lambda})$. 
Set $\myP(\lambda) := \myP(\mathcal{D}_{\lambda}) = \lambda(\upsilon\myP(\mathcal{C}))$, $\myS(\lambda) := \myS(\mathcal{D}_{\lambda}) = \lambda(\upsilon\myS(\mathcal{C}))$, and $\myA(\lambda) := \myA(\mathcal{D}_{\lambda}) = \lambda^{2}(\upsilon^{2}\myA(\mathcal{C}))$. 
Then $\myA'(\lambda) = 2\lambda(\upsilon^{2}\myA(\mathcal{C})) = 2\lambda(\upsilon\myS(\mathcal{C}))$. 
So, 
\begin{equation}
\myA'(\lambda) = \frac{d\myA}{d\lambda} = 2\myS(\lambda) = \myP(\lambda).
\end{equation}
That is, when we use the ratio of area--to--semiperimeter as the indexing parameter for a family of similar shapes, then perimeter is the derivative of area with respect to this parameter.  
Here we adopt the nomenclature of \cite{BYU} and refer to this parameter as the {\em Tong inradius} for such a family of similar shapes. 

In this paper, we use the perspective of \cite{Tong} to develop the apparently new concept of a unit shape within any given family of similar shapes. 
We will see that the relationship of this unit shape to its other similar shapes strongly analogizes the relationship of the unit circle to all circles. 
As examples, we will characterize unit right triangles, unit triangles in general, unit rectangles, unit rhombi, unit parallelograms, and unit ellipses.  
Further, we consider extremal problems relating to unit shapes, such as: Which unit triangle encloses the minimum area? 
We demonstrate that for many such questions, there are equivalent isoperimetric problems. 
We propose what we believe is naturally unital about each of the unit shapes considered here and remark on how such notions might contribute to the $\pi$ vs.\ $\tau$ debate.  
Throughout we propose many related calculus problems, some of which are open.  
We formulate all these ideas using the language of first-, second-, and third-semester Calculus \& Analytic Geometry together with some elementary linear algebra.

{\bf \S 2\ \ Shapes, families of shapes, and similarity.} 
We begin by setting up some preliminary geometric notions framed in the language of elementary calculus and linear algebra. 
Our work will mostly take place in the Euclidean plane, which is here denoted $\mathbb{R}^{2}$. 
A {\em curve} is some continuous function $f: I \longrightarrow \mathbb{R}^{2}$ from an interval $I$ on the real line to the Euclidean plane with coordinate functions, perhaps designated as  $x_{f}$ and $y_{f}$, such that $f(t) = (x_{f}(t),y_{f}(t))$. 
When such a curve is continuously differentiable on $[a,b]$ (taking a right-derivative at $a$ and a left-derivative at $b$) such that the derivatives $x_{f}'(t)$ and $y_{f}'(t)$ do not simultaneously vanish, then $\int_{a}^{b}\sqrt{x_{f}'(t)^{2}+y_{f}'(t)^{2}}\, dt$ is some positive real number; we call $f$ {\em nicely rectifiable}. 
The curve $f$ is {\em piecewise nice} if it is continuous on $[a,b]$ and there is some partition $a=a_{0} < a_{1} < \cdots < a_{n-1} < a_{n} = b$ of the interval such that the restriction of $f$ to each subinterval $[a_{i-1},a_{i}]$ (i.e.\ each ``piece'' of $f$) is nicely rectifiable. 
For us, a {\em shape} $\mathcal{C}$ will be a piecewise nice simple closed curve in the Euclidean plane, notated as $\mathcal{C}: [a,b] \longrightarrow \mathbb{R}^{2}$ with $t \stackrel{\mathcal{C}}{\mapsto} (x_{\mathcal{C}}(t),y_{\mathcal{C}}(t))$ and $\mathcal{C}(a) = \mathcal{C}(b)$. 
Most often we identify the shape (or curve) $\mathcal{C}$ with its image $\mathcal{C}([a,b])$ in the plane. 

A {\em rigid motion} $M: \mathbb{R}^{2} \longrightarrow \mathbb{R}^{2}$ of the Euclidean plane $\mathbb{R}^{2}$ is the composition of a rotation $R$, (possibly) followed by a reflection $S$, and followed by a translation $T$. 
Two shapes $\mathcal{C}$ and $\mathcal{D}$ are {\em congruent}, written $\mathcal{C} \cong \mathcal{D}$, if $\mathcal{D} = M(\mathcal{C})$ for some rigid motion $M$. 
A rotation $R$ of the plane through an angle $\theta$ can be expressed via matrix multiplication: 
\[\left(\begin{array}{c}x\\ y\end{array}\right) \stackrel{R}{\longmapsto} \left(\begin{array}{cc}\cos \theta & -\sin \theta\\ \sin \theta & \cos \theta\end{array}\right)\left(\begin{array}{c}x\\ y\end{array}\right),\] 
where a rotation through $\theta=0$ corresponds to the identity transformation of the plane. 
A reflection $S$ of the plane across a line through the origin in the direction of some unit vector $\left(\begin{array}{c}a\\ b\end{array}\right)$ can be expressed thusly: 
\[\left(\begin{array}{c}x\\ y\end{array}\right) \stackrel{S}{\longmapsto} \left(\begin{array}{cc}a^{2}-b^{2} & 2ab\\ 2ab & b^{2}-a^{2}\end{array}\right)\left(\begin{array}{c}x\\ y\end{array}\right).\] 
A translation $T$ of the plane by a vector $\left(\begin{array}{c}h\\ k\end{array}\right)$ is given by 
\[\left(\begin{array}{c}x\\ y\end{array}\right) \stackrel{T}{\longmapsto} \left(\begin{array}{c}x\\ y\end{array}\right) + \left(\begin{array}{c}h\\ k\end{array}\right) = \left(\begin{array}{c}x+h\\ y+k\end{array}\right).\] 
For clarity, we sometimes write $T = T_{(h,k)}$.  
The rigid motion $M$ may or may not include a reflection $S$; we use a binary exponent $\epsilon_{\mbox{\tiny $M$}} = 1$ to indicate the inclusion of $S$ when we write $M$ as a composite function and $\epsilon_{\mbox{\tiny $M$}} = 0$ to indicate the exclusion of $S$. 
So, $M = T_{(h,k)} \circ S^{\epsilon_{\mbox{\tiny $M$}}} \circ R$. 
We leave it as an exercise for the reader to confirm that if $M$ and $N$ are rigid motions of the plane, then so is $M^{-1} \circ N$. 

Let $\lambda$ be a positive real number.  
The {\em planar scaling function} $L_{\lambda}: \mathbb{R}^{2} \longrightarrow \mathbb{R}^{2}$ is given by $L_{\lambda}(x,y) := (\lambda x, \lambda y)$. (The notation ``$L$,'' which abbreviates the word ``length-modifying,'' is our mnemonic for remembering that $L_{\lambda}$ is scaling function.) 
Observe that for a rigid motion $M = T_{(h,k)} \circ S^{\epsilon_{\mbox{\tiny $M$}}} \circ R$, we have $L_{\lambda} \circ M = M' \circ L_{\lambda}$, where $M' = T_{(\lambda h,\lambda k)} \circ S^{\epsilon_{\mbox{\tiny $M$}}} \circ R$. 
We write $\mathcal{C} \sim_{\lambda} \mathcal{D}$ and say shape $\mathcal{D}$ is $\lambda$-{\em similar} to shape $\mathcal{C}$ if $\mathcal{D} = (L_{\lambda} \circ M)(\mathcal{C})$ for some rigid motion $M$. 
Shapes $\mathcal{C}$ and $\mathcal{D}$ are {\em similar} if $\mathcal{C} \sim_{\kappa} \mathcal{D}$ for some positive real number $\kappa$. 

The famous Jordan Curve Theorem asserts that the complement of a given shape $\mathcal{C}$ in the Euclidean plane consists of two disjoint open subsets, one of which is bounded (see e.g.\ \cite{Hales}). 
We let $\myR(\mathcal{C})$ denote the latter region together with $\mathcal{C}$ and refer to $\myR(\mathcal{C})$ as the {\em region of the plane enclosed by the shape} $\mathcal{C}$, which therefore has $\mathcal{C}$ has its boundary. 
Now, Green's Theorem affords the following identity:  
\[\iint_{\mysmallerR(\mathcal{C})}\, 1\, dx\, dy = \frac{1}{2}\ointctrclockwise_{\mathcal{C}}(x\, dy -y\, dx),\]
where the latter is a line integral computed via a counterclockwise parameterization of the curve $\mathcal{C}$. 
We therefore let $\myA(\mathcal{C}) := \iint_{\mysmallerR(\mathcal{C})}\, 1\, dx\, dy > 0$ be the {\em area} of our region $\myR(\mathcal{C})$. 
In addition, we let $\myP(\mathcal{C}) := \sum_{i=1}^{n}\int_{a_{i-1}}^{a_{i}}\sqrt{x_{\mathcal{C}}'(t)^{2}+y_{\mathcal{C}}'(t)^{2}}\, dt > 0$ be the {\em perimeter} of $\mathcal{C}$, where the $n$ subintervals of the partition $a=a_{0} < a_{1} < \cdots a_{n-1} < a_{n} = b$ of the interval $[a,b]$ correspond to the $n$ nicely rectifiable pieces of the shape $\mathcal{C}$; and we let $\myS(\mathcal{C}) := \frac{1}{2}\myP(\mathcal{C})$ be the {\em semiperimeter} of $\mathcal{C}$. 
Within the context of unit shapes (a concept we develop in the next section), the authors have come to view the semiperimeter of a shape as a measure more mathematically intrinsic than the shape's perimeter, so this quantity will be prominently featured in what follows. 

Next is a standard result about the behavior of area and semiperimeter under scaling. 
Our first calculus problems ask the reader to confirm this result in a specialized Calculus I context and then in the more general Calculus III setting.  

\noindent 
{\bf Proposition  2.1}\ \ {\sl Suppose a shape $\mathcal{D}$ is $\lambda$-similar to a shape $\mathcal{C}$. Then} 
$\myA(\mathcal{D}) = \lambda^{2}\myA(\mathcal{C})$ {\sl and} $\myS(\mathcal{D}) =  \lambda \myS(\mathcal{C})$.

\noindent 
{\bf Problems 2.2}\ \ (a) In a Calculus I context, suppose $\mathcal{C}$ is the boundary of the planar region bounded by $x=a$, $y=0$, $x=b$, and $y=f(x)$, where $f$ is a positively-valued differentiable function on the open interval $(a,b)$ that is right-differentiable at $x=a$ and left-differentiable at $x=b$. 
So, $\displaystyle \myA(\mathcal{C}) = \int_{a}^{b} f(x)dx$ and $\displaystyle \myS(\mathcal{C}) = f(a) + f(b) + (b-a) + \int_{a}^{b}\sqrt{1 + f'(x)^{2}\, }dx$. 
For a positive real number $\lambda$, let $\mathcal{D}$ be the image of $\mathcal{C}$ when the planar scaling function $L_{\lambda}$ is applied. 
Find integration-based formulas for $\myA(\mathcal{D})$ and $\myS(\mathcal{D})$, and then verify that $\myA(\mathcal{D}) = \lambda^{2}\myA(\mathcal{C})$ and $\myS(\mathcal{D}) =  \lambda \myS(\mathcal{C})$.
(b) Now within a Calculus III setting, say we have $\mathcal{D} = L_{\lambda} \circ M \circ \mathcal{C}: [a,b] \longrightarrow \mathbb{R}^{2}$ for a rigid motion $M$ and planar scaling function $L_{\lambda}$. 
Confirm that the Jacobian of the transformation $L_{\lambda} \circ M$ that sends region $\myR(\mathcal{C})$ to $\myR(\mathcal{D})$ is $\lambda^{2}$. 
Why does it follow that $\myA(\mathcal{D}) = \lambda^{2}\myA(\mathcal{C})$? 
Interior to a subinterval $[a_{i-1},a_{i}]$ of $[a,b]$, use the Chain Rule to check that $x'_{\mathcal{D}}(t) = \frac{d}{dt}\left(\rule[-1.5mm]{0mm}{5mm}L_{\lambda} \circ M \circ x_{\mathcal{C}}(t)\right) = \lambda x'_{\mathcal{C}}(t)$ and similarly that $y'_{\mathcal{D}} =  \lambda y'_{\mathcal{C}}$. 
Why does it follow that $\myS(\mathcal{D}) = \lambda \myS(\mathcal{C})$?


Given an indexing set $\mathcal{I}$ (such as the set $\mathbb{R}_{+}$ of positive real numbers), an $\mathcal{I}$-{\em family of shapes} $\{\mathcal{F}_{\lambda}\}_{\lambda \in \mathcal{I}}$ is a collection for which $\mathcal{F}_{\kappa}$ is a shape for each $\lambda \in \mathcal{I}$. 
A {\em properly similar indexing} of an  $\mathbb{R}_{+}$-family of similar shapes $\{\mathcal{F}_{\lambda}\}_{\lambda \in \mathbb{R}_{+}}$ has the property that for any $\kappa, \mu \in \mathbb{R}_{+}$, then $\mathcal{F}_{\kappa} \sim_{\mu/\kappa} \mathcal{F}_{\mu}$, in which case we say that $\{\mathcal{F}_{\lambda}\}_{\lambda \in \mathbb{R}_{+}}$ is a {\em properly similar family of shapes}. 

\noindent 
{\bf Example 2.3}\ \ For any positive real number $\lambda$, let $\mathcal{C}_{\lambda}$ be the right isosceles triangle with vertices at the origin, $(\lambda,0)$, and $(0,\lambda)$. 
We leave it as an exercise for the reader to check that $\{\mathcal{C}_{\lambda}\}_{\lambda \in \mathbb{R}_{+}}$ is a properly similar family of shapes.

As an example of the usefulness of the preceding language, we present a generalization of the Pythagorean Theorem. 
This is a version of the so-called ``Blob Pythagorean Theorem'' and was essentially known to Euclid. 
In what follows, positive real numbers $a$, $b$, and $c$ form a {\em right triangle triple} $(a,b,c)$ if there exists a right triangle with legs of lengths $a$ and $b$ and hypotenuse with length $c$. 
Note that the theorem statement does not refer to perimeter or semiperimeter, whose definitions require the Pythagorean Theorem. 

\noindent 
{\bf Theorem 2.4 (The Blob Pythagorean Theorem)}\ \ {\sl Suppose $a, b, c \in \mathbb{R}_{+}$.  Then the following are equivalent:}\\
{\sl (1) There exists a properly similar family of shapes $\{\mathcal{C}_{\lambda}\}_{\lambda \in \mathbb{R}_{+}}$ such that} $\myA(\mathcal{C}_{a}) + \myA(\mathcal{C}_{b}) = \myA(\mathcal{C}_{c})$.\\
{\sl (2) $a^{2} + b^{2} = c^{2}$.}\\
{\sl (3) For every properly similar family of shapes $\{\mathcal{C}_{\lambda}\}_{\lambda \in \mathbb{R}_{+}}$, we have} $\myA(\mathcal{C}_{a}) + \myA(\mathcal{C}_{b}) = \myA(\mathcal{C}_{c})$.\\ 
{\sl (4) $(a,b,c)$ is a right triangle triple.} 

{\em ``Proof.''} Barry Mazur renders this theorem as a fable in the very entertaining video \cite{Mazur}. 
We urge the reader to consult that video for a wonderfully clever proof that is based on what is sometimes called ``Einstein's proof by dissection without rearrangement.''\hfill\QED

{\bf \S 3\ \ Calculus-friendly indexing and unit shapes.} 
Given a shape $\mathcal{C}$, the simplest way to construct a properly similar family of shapes is to declare, for each $\lambda \in \mathbb{R}_{+}$, that $\mathcal{C}_{\lambda} := L_{\lambda}(\mathcal{C})$. 
We will refer to $\{\mathcal{C}_{\lambda}\}_{\lambda \in \mathbb{R}_{+}}$ as the {\em standard indexing} of the $\mathbb{R}_{+}$-family of shapes similar to $\mathcal{C}$; clearly the standard indexing is properly similar. 
We say a properly similar indexing $\{\mathcal{C}_{\lambda}\}_{\lambda \in \mathbb{R}_{+}}$ is {\em calculus-friendly} if, when we let $\myA(\lambda) := \myA(\mathcal{C}_{\lambda})$ and $\myS(\lambda) := \myS(\mathcal{C}_{\lambda})$, then $\myA'(\lambda) = 2\myS(\lambda)$. 
This definition is motivated by fact (1) from \S 1, which is, in effect, the observation that the indexing of origin-centered circles by their radii is calculus-friendly.  
The next theorem, whose calculus content is mostly contained in its invocations of Proposition 2.1, lays the groundwork for our notion of unit shapes. 

\noindent
{\bf Theorem 3.1}\ \ {\sl Fix a shape $\mathcal{C}$ and its (properly similar) standard indexing $\{\mathcal{C}_{\lambda}\}_{\lambda \in \mathbb{R}_{+}}$.  Let} $\ddot{\upsilon} := \myS(\mathcal{C})/\myA(\mathcal{C})$. 
{\sl (1) For each positive real number $\lambda$, let $\mathcal{D}_{\lambda}$ be the shape $\mathcal{C}_{\lambda \ddot{\upsilon}}$. Then $\{\mathcal{D}_{\lambda}\}_{\lambda \in \mathbb{R}_{+}}$ is a calculus-friendly properly similar indexing of this family of shapes. Moreover, $\mathcal{D}_{1} = \mathcal{C}_{\ddot{\upsilon}}$ is the unique shape $\mathcal{X}$ within the family $\{\mathcal{D}_{\lambda}\}_{\lambda \in \mathbb{R}_{+}}$ that satisfies} $\myA(\mathcal{X}) = \myS(\mathcal{X})$, {\sl and $\{\mathcal{D}_{\lambda}\}_{\lambda \in \mathbb{R}_{+}}$ is the standard indexing of the family of shapes similar to $\mathcal{D}_{1}$. 
(2) Now suppose $\{\mathcal{B}_{\lambda}\}_{\lambda \in \mathbb{R}_{+}}$ is another calculus-friendly properly similar indexing of a family of shapes such that for some positive real number $\kappa$ we have $\mathcal{B}_{\kappa} \cong \mathcal{C}$. Then $\mathcal{B}_{\lambda} \cong \mathcal{D}_{\lambda}$ for each positive real number $\lambda$, and} $\lambda = \myA(\mathcal{B}_{\lambda})/\myS(\mathcal{B}_{\lambda})$. 
{\sl (3) If $\mathcal{U}$ is any shape similar to $\mathcal{C}$ such that} $\myA(\mathcal{U}) = \myS(\mathcal{U})$, {\sl then $\mathcal{U} \cong \mathcal{D}_{1} = \mathcal{C}_{\ddot{\upsilon}}$.}

Before we present the proof, we offer two remarks. 
First, we comment on some notation introduced in part {\sl (1)} of the theorem statement. 
Since the Greek letter upsilon ``$\upsilon$'' seems appropriate within our context but is, visually, a rather undistinguished symbol, we have utilized a rare umlaut variation of this letter: ``$\, \ddot{\upsilon}\, $.'' 
This notation also announces the happy fact that $\mathcal{C}_{\ddot{\upsilon}}$ is what we will shortly refer to as a ``unit shape.''  
Second, note that part {\sl (2)} of the theorem (re)asserts that the only calculus-friendly indexing parameter for a properly similar family of shapes is the Tong inradius. 

{\em Proof of Theorem 3.1.} For {\sl (1)}, consider two positive real numbers $\kappa$ and $\lambda$. 
Then \[L_{\lambda/\kappa}(\mathcal{D}_{\kappa}) = L_{\lambda/\kappa}(\mathcal{C}_{\kappa \ddot{\upsilon}}) = (L_{\lambda/\kappa} \circ L_{\kappa \ddot{\upsilon}})(\mathcal{C}) = L_{\lambda \ddot{\upsilon}}(\mathcal{C}) = \mathcal{C}_{\lambda \ddot{\upsilon}} = \mathcal{D}_{\lambda},\]
and hence $\mathcal{D}_{\kappa} \sim_{\lambda/\kappa} \mathcal{D}_{\lambda}$. 
So, $\{\mathcal{D}_{\lambda}\}_{\lambda \in \mathbb{R}_{+}}$ is a properly indexed family of similar shapes. 
Now $\mathcal{D}_{1} \sim_{\lambda} \mathcal{D}_{\lambda}$, and hence $\myA(\mathcal{D}_{\lambda}) = \lambda^{2}\myA(\mathcal{D}_{1})$ and $\myS(\mathcal{D}_{\lambda}) = \lambda\myS(\mathcal{D}_{1})$ by Proposition 2.1. 
Also, note that $\mathcal{D}_{1} = \mathcal{C}_{\ddot{\upsilon}} \sim_{1/\ddot{\upsilon}} \mathcal{C}$. 
So: 
$\myA(\mathcal{D}_{1})
 = {\ddot{\upsilon}^{2}}\myA(\mathcal{C})
 =  \frac{\myS(\mathcal{C})^{2}}{\myA(\mathcal{C})^{2}}\myA(\mathcal{C})
 = \frac{\myS(\mathcal{C})}{\myA(\mathcal{C})}\myS(\mathcal{C})
 = {\ddot{\upsilon}}\myS(\mathcal{C}) 
 = \myS(\mathcal{D}_{1})$. 
Therefore $\frac{d}{d\lambda}\left(\rule[-1.5mm]{0mm}{5mm}\myA(\mathcal{D}_{\lambda})\right) = \frac{d}{d\lambda}\left(\rule[-1.5mm]{0mm}{5mm}\lambda^{2}\myA(\mathcal{D}_{1})\right) = 2\lambda\myA(\mathcal{D}_{1}) = 2\lambda\myS(\mathcal{D}_{1}) = 2\myS(\mathcal{D}_{\lambda})$, so the properly similar indexing $\{\mathcal{D}_{\lambda}\}_{\lambda \in \mathbb{R}_{+}}$ is calculus-friendly. 
Finally, suppose $\myA(\mathcal{D}_{\kappa}) = \myS(\mathcal{D}_{\kappa})$ for some positive real number $\kappa$. 
Since $\mathcal{D}_{\kappa} \sim_{1/\kappa} \mathcal{D}_{1}$, then $\myA(\mathcal{D}_{\kappa}) = \kappa^{2}\myA(\mathcal{D}_{1}) = \kappa^{2}\myS(\mathcal{D}_{1}) = \kappa\myS(\mathcal{D}_{\kappa}) = \kappa\myA(\mathcal{D}_{\kappa})$, so we must have $\kappa=1$. 

For {\sl (2)}, we have $\mathcal{B}_{\kappa} \sim_{\lambda/\kappa} \mathcal{B}_{\lambda}$ for any positive real number $\lambda$, so $\mathcal{B}_{\lambda} = (L_{\lambda/\kappa} \circ M)(\mathcal{B}_{\kappa})$ for some rigid motion $M$. 
Also, since $\mathcal{B}_{\kappa} \cong \mathcal{C}$ then we have $\mathcal{B}_{\kappa} = N(\mathcal{C})$ for some rigid motion $N$, and hence $\mathcal{B}_{\lambda} = (L_{\lambda/\kappa} \circ M \circ N)(\mathcal{C})$. 
As noted in \S 2, there are rigid motions $M'$ and $N'$ such that $L_{\lambda/\kappa} \circ M \circ N = M' \circ N' \circ L_{\lambda/\kappa}$. 
Then $\mathcal{B}_{\lambda} = (L_{\lambda/\kappa} \circ M \circ N)(\mathcal{C}) = (M' \circ N')\left(\rule[-1.5mm]{0mm}{5mm}L_{\lambda/\kappa}(\mathcal{C})\right) = (M' \circ N')(\mathcal{C}_{\lambda/\kappa}) = (M' \circ N')(\mathcal{D}_{\lambda/(\ddot{\upsilon} \kappa)})$, and so $\mathcal{B}_{\lambda} \cong \mathcal{D}_{\lambda/(\ddot{\upsilon} \kappa)}$. 
Next, we observe that $\myA(\mathcal{B}_{\lambda}) = \lambda^{2}\myA(\mathcal{B}_{1})$ and $\myS(\mathcal{B}_{\lambda}) = \lambda\myS(\mathcal{B}_{1})$, since $\mathcal{B}_{1} \sim_{\lambda} \mathcal{B}_{\lambda}$. 
From the fact that the properly similar indexing $\{\mathcal{B}_{\lambda}\}_{\lambda \in \mathbb{R}_{+}}$ is calculus-friendly, we obtain \[2\lambda\myA(\mathcal{B}_{1}) = \frac{d}{d\lambda}\left(\rule[-1.5mm]{0mm}{5mm}\lambda^{2}\myA(\mathcal{B}_{1})\right)  = \frac{d}{d\lambda}\left(\rule[-1.5mm]{0mm}{5mm}\myA(\mathcal{B}_{\lambda})\right) = 2\myS(\mathcal{B}_{\lambda}) = 2\lambda\myS(\mathcal{B}_{1}),\] 
and we can conclude that $\myA(\mathcal{B}_{1}) = \myS(\mathcal{B}_{1})$. 
But, since the indexing is proper, we have $\mathcal{B}_{\kappa} \sim_{1/\kappa} \mathcal{B}_{1}$, so by Proposition 2.1, $\myA(\mathcal{B}_{1}) = \frac{1}{\kappa^{2}}\myA(\mathcal{B}_{k}) = \frac{1}{\kappa^{2}}\myA(\mathcal{C})$ and $\myS(\mathcal{B}_{1}) = \frac{1}{\kappa}\myS(\mathcal{B}_{k}) = \frac{1}{\kappa}\myS(\mathcal{C})$. 
Therefore $\frac{1}{\kappa} = \frac{\mysmallerS(\mathcal{C})}{\mysmallerA(\mathcal{C})} = \ddot{\upsilon}$. 
Thus, $\mathcal{B}_{\lambda} \cong \mathcal{D}_{\lambda/(\ddot{\upsilon} \kappa)} = \mathcal{D}_{\lambda}$. 
Of course, $\myA(\mathcal{B}_{\lambda})/\myS(\mathcal{B}_{\lambda}) = \lambda^{2}\myA(\mathcal{B}_{1})/(\lambda \myS(\mathcal{B}_{1})) = \lambda$. 

For {\sl (3)}, let $\{\mathcal{B}_{\lambda}\}_{\lambda \in \mathbb{R}_{+}}$ be the standard indexing of the shape $\mathcal{U}$, where we identify $\mathcal{B}_{1}$ as $\mathcal{U}$. 
Now, by Proposition 2.1, $\myA(\mathcal{B}_{\lambda}) = \lambda^{2}\myA(\mathcal{U})$ and $\myS(\mathcal{B}_{\lambda}) = \lambda\myS(\mathcal{U})$, so $\frac{d}{d\lambda}\left(\rule[-1.5mm]{0mm}{5mm}\myA(\mathcal{B}_{\lambda})\right) = \frac{d}{d\lambda}\left(\rule[-1.5mm]{0mm}{5mm}\lambda^{2}\myA(\mathcal{U})\right) = 2\lambda\myA(\mathcal{U}) = 2\lambda\myS(\mathcal{U})$. 
We conclude that the standard (properly similar) indexing $\{\mathcal{B}_{\lambda}\}_{\lambda \in \mathbb{R}_{+}}$ is calculus-friendly. 
By hypothesis, there exists a positive real number $\kappa$ and a rigid motion $M$ such that $\mathcal{C} = (L_{\kappa} \circ M)(\mathcal{U})$. 
So for some rigid motion $M'$ we get $\mathcal{C} = (L_{\kappa} \circ M)(\mathcal{U}) = (M' \circ L_{\kappa})(\mathcal{U}) = M'(\mathcal{B}_{\kappa})$. 
That is, $\mathcal{B}_{\kappa} \cong \mathcal{C}$. 
So we meet the criteria for part {\sl (2)} of the theorem statement, and hence $\mathcal{U} = \mathcal{B}_{1} \cong \mathcal{D}_{1} = \mathcal{C}_{\ddot{\upsilon}}$.\hfill\QED

In view of the preceding theorem, we make the following definitions and observations. 

\noindent 
{\bf Definitions/Observations 3.2}\ \ A {\em unit shape} is a shape $\mathcal{U}$ for which $\myA(\mathcal{U}) = \myS(\mathcal{U})$.  
By Theorem 3.1, any two unit shapes which are similar must be congruent.  
For any given shape $\mathcal{C}$, a $\mathcal{C}$-{\em similar unit shape} is any unit shape that is similar to $\mathcal{C}$. 
By Theorem 3.1, any $\mathcal{C}$-similar unit shape is congruent to $\canonicalU(\mathcal{C}) := L_{\ddot{\upsilon}}(\mathcal{C})$, where $\ddot{\upsilon} := \myS(\mathcal{C})/\myA(\mathcal{C})$, and the standard indexing $\{\canonicalU(\mathcal{C})_{_\lambda}\}_{\lambda \in \mathbb{R}_{+}}$ of the family of shapes similar to $\canonicalU(\mathcal{C})$ is calculus-friendly. 
For these reasons, we call $\canonicalU(\mathcal{C})$ the {\em canonical $\mathcal{C}$-similar unit shape}. 
Finally, we declare that $\Pi_{\mathscr{U}\!(\mathcal{C})} := \myA(\canonicalU(\mathcal{C})) = \myS(\canonicalU(\mathcal{C}))$ and call this quantity the {\em fundamental measure for shapes similar to} $\mathcal{C}$. 
If $\mathcal{C}$ itself is a unit shape, we sometimes use the more efficient notation $\Pi_{\mathcal{C}}$ as a substitute for $\Pi_{\mathscr{U}\!(\mathcal{C})}$ and call this the {\em fundamental measure of the unit shape} $\mathcal{C}$. 
By Proposition 2.1, $\myA(\canonicalU(\mathcal{C})_{_\lambda}) = \lambda^{2}\Pi_{\mathscr{U}\!(\mathcal{C})}$ and $\myS(\canonicalU(\mathcal{C})_{_\lambda}) = \lambda\Pi_{\mathscr{U}\!(\mathcal{C})}$, and hence $\frac{d}{d\lambda}\left(\rule[-1.5mm]{0mm}{5mm}\myA(\canonicalU(\mathcal{C})_{_\lambda})\right) = 2\myS(\canonicalU(\mathcal{C})_{_\lambda})$.\hfill\QED 

The remainder of the paper is devoted to exploring this unit shape concept and will hopefully make the case that this notion yields many interesting and illustrative first-, second-, and third-semester Calculus examples and problems. 
To close this section, we invite the reader to reconsider the opening paragraph of this paper in view of unit shapes. 

\noindent 
{\bf Problem 3.3}\ \ Let $\mathcal{U}$ be a unit shape, $\{\mathcal{U}_{\lambda}\}_{\lambda \in \mathbb{R}_{+}}$ the calculus-friendly properly similar family of shapes with $\mathcal{U}_{1} = \mathcal{U}$, $\myA(\lambda)$ the area function affiliated with this family, and $\myS(\lambda)$ the related semiperimeter function. 
Generalize the argument from the first paragraph of \S 1 relating the quantities $\Delta \myA$, $d\myA$, $d\lambda$, and $\myS(\lambda)$.

{\bf \S 4\ \ A menagerie of unit shapes.} Here we introduce a number of examples of unit shapes, and, as we proceed, we encourage the reader to consider the question: What is \underline{unital} about each of these shapes?

\vspace*{0.1in}
\noindent 
\parbox[b]{5.6in}{{\bf Example 4.1: Unit right triangles}\ \ Fix an acute angle measure $\theta$, so $0 < \theta < \pi/2$. 
For any given positive real number $b$, let $\mathcal{N}^{(\theta)}_{b}$ be a right triangle with a leg of length $b$ adjacent to an acute angle of measure $\theta$ (depicted to the right), where we think of $b$ as a base measurement. 
(Here, the letter ``$\mathcal{N}$'' is indicative of the \underline{\em n}inety-degree angle possessed by each right triangle.)  
Then $\mathcal{N}^{(\theta)}_{b}$ has area $\frac{1}{2}b^{2} \tan\theta$ and semiperimeter $\frac{1}{2}b\left(\rule[-1.5mm]{0mm}{5mm}1+\tan\theta+\sec\theta\right)$. Therefore the canonical $\mathcal{N}^{(\theta)}_{b}$-similar unit triangle $\mathcal{N}^{(\theta)} := \canonicalU(\mathcal{N}^{(\theta)}_{b})$ has base length $1+\cot\theta+\csc\theta$. Moreover, the fundamental measure of $\mathcal{N}^{(\theta)}$ is $\displaystyle \Pi_{\mathcal{N}^{(\theta)}} = $ $\frac{1}{2}\left(\rule[-1.5mm]{0mm}{5mm}1+\tan\theta +\sec\theta \right)\left(\rule[-1.5mm]{0mm}{5mm}1+\cot\theta+\csc\theta\right) = \left(\rule[-1.5mm]{0mm}{5mm}1+\sec\theta\right)\left(\rule[-1.5mm]{0mm}{5mm}1+\csc\theta\right)$.} 
\setlength{\unitlength}{1cm}
\begin{picture}(0,0)
\put(-0.2,0.2){\includegraphics[scale=0.9]{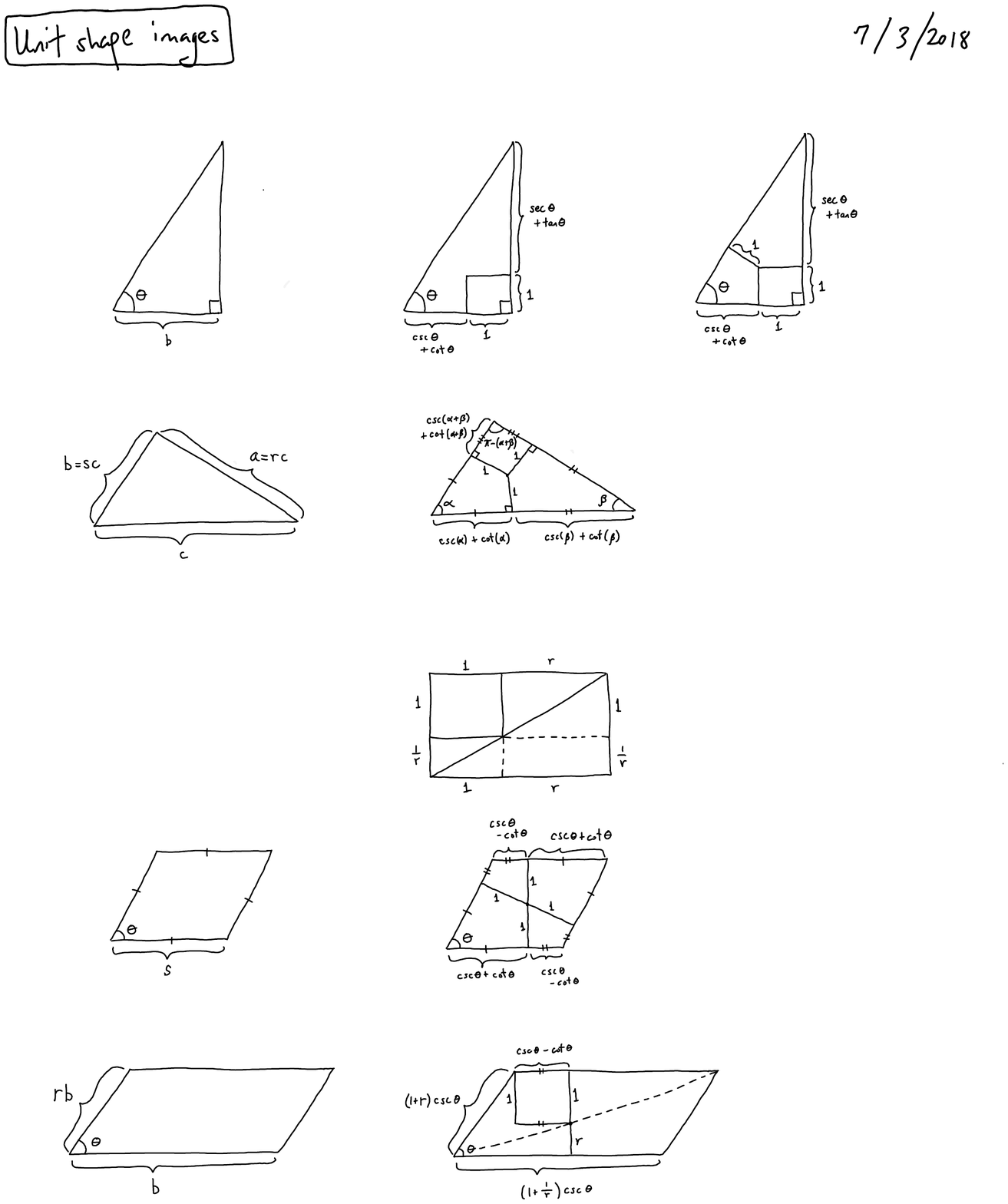}}
\end{picture}

\noindent 
{\bf Problem 4.2}\ \ Use first-semester Calculus methods to determine which acute angle measure $\theta$ minimizes the area/semiperimeter amongst all unit triangles $\{\mathcal{N}^{(\theta)}\}_{0 < \theta < \pi/2}$. 
What is this minimum fundamental measure? 
Describe the associated shape using familiar language. 

\vspace*{0.1in}
\noindent 
\parbox[b]{4.75in}{{\bf Example 4.3: Unit triangles, in general}\ \ Now fix two positive real numbers $\myr$ and $\mys$ to be thought of as ratios satisfying the inequalities $\myr \leq 1$, $\mys \leq 1$, and $\myr+\mys > 1$.  
For reasons that will be clear momentarily, we call such a ratio pair $(\myr,\mys)$ {\em triangle-friendly}. 
For any given positive real number $c$, let $\mathcal{T}^{(\mynewerr,\mynewers)}_{c}$ denote a triangle with sides of lengths $a=\myr c$, $b=\mys c$, and $c$, with longest side of length $c$.  
See the depiction to the right. 
(Why is such a triangle guaranteed to exist?) 
Then $\mathcal{T}^{(\mynewerr,\mynewers)}_{c}$ has semiperimeter $\myS := \frac{c}{2}(\myr+\mys+1)$ and, by Heron's formula, area} 
\setlength{\unitlength}{1cm}
\begin{picture}(0,0)
\put(0,1){\includegraphics[scale=0.9]{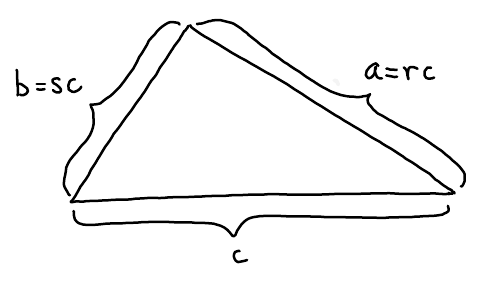}}
\end{picture} 
\[\sqrt{\myS(\myS-a)(\myS-b)(\myS-c)\, } = \frac{c^{2}}{4}\sqrt{(\myr+\mys+1)(-\myr+\mys+1)(\myr-\mys+1)(\myr+\mys-1)\, }.\]
So the canonical $\mathcal{T}^{(\mynewerr,\mynewers)}_{c}$-similar unit triangle $\mathcal{T}^{(\mynewerr,\mynewers)} := \canonicalU(\mathcal{T}^{(\mynewerr,\mynewers)}_{c})$ has $c=2\sqrt{\frac{\myr+\mys+1}{(-\myr+\mys+1)(\myr-\mys+1)(\myr+\mys-1)}\, }$  
and fundamental measure  
$\displaystyle \Pi_{\mathcal{T}^{(\mynewerr,\mynewers)}} = \frac{(\myr+\mys+1)^{3/2}}{(-\myr+\mys+1)^{1/2}(\myr-\mys+1)^{1/2}(\myr+\mys-1)^{1/2}}$.

\noindent 
{\bf Problems 4.4}\ \ (a) Conciliate our formula for the fundamental measure $\Pi_{\mathcal{T}^{(\mynewerr,\mynewers)}}$ from Example 4.3 with our formula for the fundamental measure $\Pi_{\mathcal{N}^{(\theta)}}$ from Example 4.1. 
(b) Use third-semester Calculus methods to determine which triangle-friendly pair minimizes the area/semiperimeter amongst all unit triangles $\{\mathcal{T}^{(\mynewerr,\mynewers)}\}_{\mbox{\scriptsize triangle-friendly pairs}\, (\mynewerr,\mynewers)}$. 
What is this minimum fundamental measure? 
Describe the associated shape using familiar language.

\noindent {\bf Example 4.5: Unit rectangles}\ \ Fix a positive real number $\myr$, to be thought of as a ratio. 
For any given positive real number $\ell$, let $\mathcal{R}^{(\mynewerr)}_{\ell}$ be a rectangle with length $\ell$ and height $h=\myr \ell$, as in the 

\newpage 
\noindent 
\parbox[b]{4.8in}{diagram to the right. 
Then $\mathcal{R}^{(\mynewerr)}_{\ell}$ has area $\myr \ell^{2}$ and semiperimeter $(1+\myr)\ell$. 
So the unit rectangle $\mathcal{R}^{(\mynewerr)} := \canonicalU(\mathcal{R}^{(\mynewerr)}_{\ell})$ has length $\frac{1+\mynewr}{\mynewr} = 1 + \frac{1}{\mynewr}$ and height $1+\myr$. 
The fundamental measure of $\mathcal{R}^{(\mynewerr)}$ is 
$\displaystyle \Pi_{\mathcal{R}^{(\mynewerr)}} = \frac{(1+\myr)^{2}}{\myr}$.} 
\setlength{\unitlength}{1cm}
\begin{picture}(0,0)
\put(0,-0.4){\includegraphics[scale=0.9]{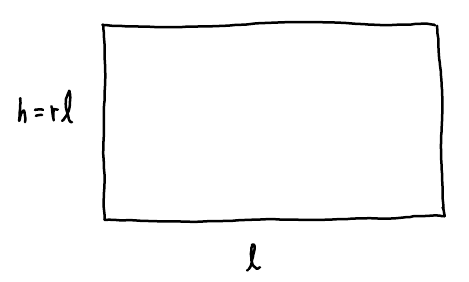}}
\end{picture} 

\vspace*{0.1in}
\noindent 
{\bf Problems 4.6}\ \ (a) What is the fundamental measure of the unit golden rectangle $\mathcal{G}$? the unit square $\mathcal{S}$? 
(b) Use first-semester Calculus methods to determine which ratio $\myr$ minimizes the area/semiperimeter amongst all unit rectangles $\{\mathcal{R}^{(\mynewerr)}\}_{\mynewerr \in \mathbb{R}_{+}}$. 
What is this minimum fundamental measure?
Describe the associated shape using familiar language. 

\vspace*{0.1in}
\noindent 
\parbox[b]{5.25in}{{\bf Example 4.7: Unit rhombi}\ \ Fix an angle measure $\theta$ with $0 < \theta < \pi$. 
For any given positive real number $s$, let $\mathcal{Q}^{(\theta)}_{s}$ be a rhombus (``e\underline{\em q}uilateral \underline{\em q}uadrilateral'') with an interior angle of measure $\theta$ and side lengths all equal to $s$. 
See the depiction to the right. 
Then $\mathcal{Q}^{(\theta)}_{s}$ has area $s^{2}\sin \theta$ and semiperimeter $2s$. 
So, the unit rhombus $\mathcal{Q}^{(\theta)} := \canonicalU(\mathcal{Q}^{(\theta)}_{s})$ has side length equal to $2\csc \theta$. 
The fundamental measure of $\mathcal{Q}^{(\theta)}$ is $\displaystyle \Pi_{\mathcal{Q}^{(\theta)}} = 4\csc \theta$.} 
\setlength{\unitlength}{1cm}
\begin{picture}(0,0)
\put(0,0.4){\includegraphics[scale=0.9]{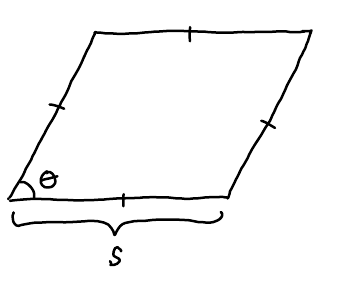}}
\end{picture} 

\noindent 
{\bf Problems 4.8}\ \ (a) Without using calculus, determine which angle measure $\theta$ minimizes the area/semiperimeter amongst all unit rhombi $\{\mathcal{Q}^{(\theta)}\}_{0 < \theta < \pi}$. 
What is this minimum fundamental measure? 
Describe the associated shape using familiar language. 
(b) Let $h(\theta)$ be the length of a shortest diagonal of the unit rhombus $\mathcal{Q}^{(\theta)}$. Determine, if possible, the greatest lower bound $l$ and the least upper bound $u$ for the set $\{h(\theta)\}_{0 < \theta < \pi}$. Is there a unit rhombus which has a shortest diagonal of length $l$? of length $u$?  

\vspace*{0.1in}
\noindent 
\parbox[b]{4.4in}{{\bf Example 4.9: Unit parallelograms}\ \ Now fix a positive real number $\myr$, to be thought of as a ratio, and an angle measure $\theta$ with $0 < \theta < \pi$. 
For any given positive real number $b$, let $\mathcal{P}^{(\theta,\mynewerr)}_{b}$ be a parallelogram with an interior angle of measure $\theta$ adjacent to two parallelogram sides of lengths $b$ and $\myr b$ respectively, where we think of $b$ as the base length measure. 
See the depiction to the right.} 
\setlength{\unitlength}{1cm}
\begin{picture}(0,0)
\put(0,0.3){\includegraphics[scale=0.9]{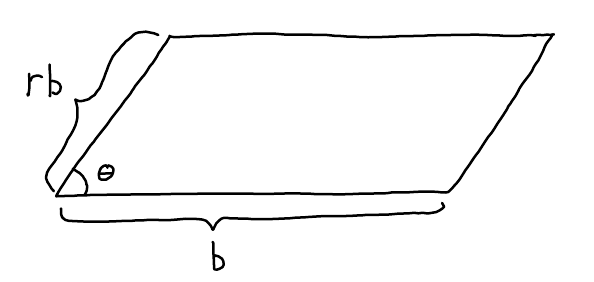}}
\end{picture}\\ 
This parallelogram has area $b^{2}\myr \sin \theta$ and semiperimeter $b(1+\myr)$. 
So, the unit parallelogram $\mathcal{P}^{(\theta,\mynewerr)} := \canonicalU(\mathcal{P}^{(\theta,\mynewerr)}_{b})$ has base length $\displaystyle \frac{1+\myr}{\, \myr \sin \theta\, }$. 
The fundamental measure of $\mathcal{P}^{(\theta,\mynewerr)}$ is  
$\displaystyle \Pi_{\mathcal{P}^{(\theta,\mynewerr)}} = \frac{(1+\myr)^{2}}{\, \myr \sin \theta\, }$.

\noindent 
{\bf Problems 4.10}\ \ (a) Conciliate our formula for the fundamental measure $\Pi_{\mathcal{P}^{(\theta,\mynewerr)}}$ from Example 4.9 with our formulas for the fundamental measures $\Pi_{\mathcal{R}^{(\mynewerr)}}$ and $\Pi_{\mathcal{D}^{(\theta)}}$ from Examples 4.5 and 4.7. 
(b) Use third-semester Calculus methods to determine which angle measure $\theta$ and ratio $\myr$ minimize the area/semiperimeter amongst all unit parallelograms $\{\mathcal{P}^{(\theta,\mynewerr)}\}_{0 < \theta < \pi, \mynewerr \in \mathbb{R}_{+}}$. 
What is this minimum fundamental measure? 
Describe the associated shape using familiar language. 

\vspace*{0.05in}
\noindent 
\parbox[b]{4.7in}{{\bf Example 4.11: Unit ellipses}\ \ Fix a positive real number $\myr$ with $0 < \myr < 1$, to be thought of as the semi-minor--to--semi-major axis ratio for an ellipse centered at the origin with foci on the $x$-axis. 
For any given positive real number $a$, let $\mathcal{E}^{(\mynewerr)}_{a}$ be the ellipse defined by the equation $\displaystyle \frac{x^{2}}{a^{2}}+\frac{y^{2}}{b^{2}} = 1$ with $b := \myr a$, as depicted to the right. 
For a natural para-}  
\setlength{\unitlength}{1cm}
\begin{picture}(0,0)
\put(0,-0.15){\includegraphics[scale=0.9]{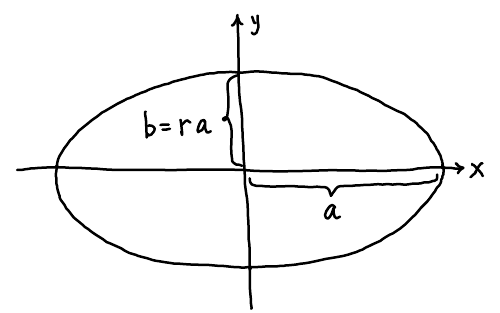}}
\end{picture}

\newpage
\noindent 
meterization of $\mathcal{E}^{(\mynewerr)}_{a}$, we use $t \mapsto \left(\rule[-1.5mm]{0mm}{5mm}x(t) := a \cos(t), y(t) := b \sin(t)\right)$ for $0 \leq t \leq 2\pi$. Then $\myA(\mathcal{E}^{(\mynewerr)}_{a}) = \pi a b = \pi \myr a^{2}$ and
$\displaystyle \myS(\mathcal{E}^{(\mynewerr)}_{a}) = \int_{0}^{\pi}\sqrt{a^{2}\sin(t)^{2}+b^{2}\cos(t)^{2}\, }dt = a\int_{0}^{\pi}\sqrt{\sin(t)^{2}+\myr^{2}\cos(t)^{2}\, }dt$. 
So the unit ellipse $\mathcal{E}^{(\mynewerr)} = \canonicalU(\mathcal{E}^{(\mynewerr)}_{a})$ has semi-major axis 
$\displaystyle \rule[-5mm]{0mm}{13mm}\frac{\int_{0}^{\pi}\sqrt{1+(\myr^{2}-1)\cos(t)^{2}\, }dt}{\pi \myr}$. 
The fundamental measure of $\mathcal{E}^{(\mynewerr)}$ is 
$\displaystyle \Pi_{\mathcal{E}^{(\mynewerr)}} = \frac{\left(\int_{0}^{\pi}\sqrt{1+(\myr^{2}-1)\cos(t)^{2}\, }dt\right)^{2}}{\pi \myr}$.

\noindent 
{\bf Problems 4.12}\ \ (a) Let $a: (0,1) \longrightarrow \mathbb{R}$ be given by $a(\myr) := \frac{1}{\pi}\int_{0}^{\pi}\sqrt{1+(\mynewerr^{2}-1)\cos(t)^{2}\, }dt$, i.e.\ the semi-minor axis of the unit ellipse $ \mathcal{E}^{(\mynewerr)}$. Is $a$ differentiable? On what subintervals is $a$ increasing? decreasing? What is the exact range of $a$? {\sc hint:} The function $a$ is bounded below by a positive (and somewhat curious) constant. 
(b) Let $f: (0,1) \longrightarrow \mathbb{R}$ be given by $f(\myr) := \Pi_{\mathcal{E}^{(\mynewerr)}} = \frac{1}{\pi \mynewerr}\left(\int_{0}^{\pi}\sqrt{1+(\mynewerr^{2}-1)\cos(t)^{2}\, }dt\right)^{2}$.  Is $f$ differentiable? On what subintervals is $f$ increasing? decreasing? What is the exact range of $f$?

{\bf \S 5\ \ Understanding isoperimetric problems in terms of unit shapes.} 
We view Viktor Bl{\aa}sj{\"o}'s paper \cite{Blasjo} as a definitive (and most entertaining) account of the history of the isoperimetric problem and liberally reference his paper in this section. 
This famous problem originated in antiquity roughly along the following lines: Amongst planar curves of some given fixed length that reside in a closed half-plane determined by some line in the plane and that adjoin the two endpoints of a given fixed segment on that line, which will maximize the enclosed area? 
An obvious variation on this question, expressed in the language of this paper, is: Which shape of some fixed semiperimeter will maximize the enclosed area? 
The well-known answer can be deduced from the next result, known as the Isoperimetric Inequality. 
In \cite{Blasjo}, Bl{\aa}sj{\"o} surveys many proofs and near-proofs of this result.  

\noindent 
{\bf Theorem 5.1 (The Isoperimetric Inequality)}\ \ {\sl For a given shape $\mathcal{C}$, we have} $\pi \myA(\mathcal{C}) \leq \myS(\mathcal{C})^{2}$, {\sl with equality if and only if $\mathcal{C}$ is a circle.}

The solution to the isoperimetric problem problem posed above now follows immediately.  
Any shape $\mathcal{C}$ with fixed semiperimeter $\myS$ has an area satisfying the inequality $\displaystyle \myA(\mathcal{C}) \leq \frac{\myS^{2}}{\pi} = \frac{r^{2}\pi^{2}}{\pi} = \pi r^{2}$, where $r := \myS/\pi$ is the radius of a circle (``halo,'' perhaps) $\mathcal{H}_{r}$ with semiperimeter $\myS$. 
Then this circle $\mathcal{H}_{r}$ is the only shape, up to congruence, that realizes the maximum possible area. 

The menagerie of \S 4 showcased unit shapes with fundamental measures (i.e.\ areas and semiperimeters) equal to $3+2\sqrt{2}$ (the unit right isosceles triangle), equal to $3\sqrt{3}$ (the unit equilateral triangle), equal to $4$ (the unit square), and greater than $\pi$ (any unit ellipse). 
Thus the natural question:\\ 
\centerline{\em Which unit shape encloses the minimum area, and what is this minimum area?} 
We will exploit the Isoperimetric Inequality to answer this question. 

\noindent 
{\bf Theorem 5.2}\ \ {\sl Let $\mathcal{H}_{1}$ denote the unit circle centered at the origin, and let $\mathcal{U}$ be any unit shape. 
Then} $\pi \leq \Pi_{\mathcal{U}}$, {\sl and} $\pi = \Pi_{\mathcal{U}}$ {\sl if and only if $\mathcal{U} \cong \mathcal{H}_{1}$.}

{\em Proof.} Now, $\pi \myA(\mathcal{U}) \leq \myS(\mathcal{U})^{2}$ by the Isoperimetric Inequality. 
But $\myS(\mathcal{U})^{2} = \myA(\mathcal{U})^{2}$, so the preceding inequality simplifies to  $\pi \leq \myA(\mathcal{U}) = \Pi_{\mathcal{U}}$. 
If $\mathcal{U} \cong \mathcal{H}_{1}$, then $\mathcal{U}$ is $1$-similar to $\mathcal{H}_{1}$, so by Proposition 2.1 we have $\myA(\mathcal{U}) = \myA(\mathcal{H}_{1}) = \pi$. 
On the other hand, suppose $\pi = \myA(\mathcal{U})$. Then $\pi \myA(\mathcal{U}) = \myA(\mathcal{U})^{2} = \myS(\mathcal{U})^{2}$, and by the Isoperimetric Inequality, $\mathcal{U}$ must be a circle. 
Since $\pi = \myA(\mathcal{U})$, then $\mathcal{U}$ has unit radius, and hence $\mathcal{U} \cong \mathcal{H}_{1}$.\hfill\QED

It is easy to see that Theorem 5.1 can be deduced from Theorem 5.2, so these two theorems are logically equivalent.  
Next, we record a simple result that connects isoperimetric problems to questions about minimizing unit shape areas. 

\noindent
{\bf Proposition 5.3}\ \ {\sl (1) Let $\rho$ be a positive real number and $\mathcal{U}$ a unit shape. 
Then} $\rho \leq \Pi_{\mathcal{U}}$ {\sl if and only if for some positive real number $\kappa$ and some shape $\mathcal{U}_{\kappa}$ that is $\kappa$-similar to $\mathcal{U}$ we have} $\rho \myA(\mathcal{U}_{\kappa}) \leq \myS(\mathcal{U}_{\kappa})^{2}$ {\sl if and only if for every positive real number $\kappa$ and any shape $\mathcal{U}_{\kappa}$ that is $\kappa$-similar to $\mathcal{U}$ we have} $\rho \myA(\mathcal{U}_{\kappa}) \leq \myS(\mathcal{U}_{\kappa})^{2}$. 
{\sl (2) Now let $\mathscr{S}$ be a collection of shapes with the property that if $\mathcal{C}$ is in $\mathscr{S}$, then for any positive real number $\lambda$ there exists $\mathcal{D}$ in $\mathscr{S}$ such that $\mathcal{C} \sim_{\lambda} \mathcal{D}$. 
Suppose $\rho$ is a positive lower bound of the fundamental measures of all unit shapes in $\mathscr{S}$. 
Then} $\rho\myA(\mathcal{C}) \leq \myS(\mathcal{C})^{2}$ {\sl for all shapes $\mathcal{C}$ in $\mathscr{S}$.} 

{\em Proof.} Let $\kappa$ be a positive real number and $\mathcal{U}_{\kappa}$ a shape that is $\kappa$-similar to $\mathcal{U}$. Then:  
$\rho \leq \Pi_{\mathcal{U}} = \myA(\mathcal{U}) 
\iff  \rho \myA(\mathcal{U}) \leq \myA(\mathcal{U})^{2} 
\iff  \rho \myA(\mathcal{U}) \leq \myS(\mathcal{U})^{2} 
\iff  \rho \lambda^{2} \myA(\mathcal{U}) \leq \lambda^{2} \myS(\mathcal{U})^{2} 
\iff \rho \myA(\mathcal{U}_{\lambda}) \leq \myS(\mathcal{U}_{\lambda})^{2}$.
This suffices to establish {\sl (1)}.  
For {\sl(2)}, suppose $\rho$ is a positive lower bound of the fundamental measures of all unit shapes in $\mathscr{S}$. 
Since any shape $\mathcal{C}$ in $\mathscr{S}$ is similar to some unit shape $\mathcal{V}$ in $\mathscr{S}$, and since $\rho \leq \Pi_{\mathcal{V}}$, then by {\sl (1)} it follows that $\rho\myA(\mathcal{C}) \leq \myS(\mathcal{C})^{2}$.\hfill\QED

As an example, consider an $m$-gon $\mathcal{M}$, where $m$ is an integer no smaller than three. 
From Bl{\aa}sj{\"o}'s account in \cite{Blasjo} of the isoperimetric problem for $m$-gons, we learn that $\rho_{m} \myA(\mathcal{M}) \leq \myS(\mathcal{M})^{2}$, where $\rho_{m} = m\tan(\pi/m)$. 
Moreover, we also learn that $\rho_{m} \myA(\mathcal{M}) = \myS(\mathcal{M})^{2}$ if and only if $\mathcal{M}$ is a regular $m$-gon. 
Let $\mathcal{U}^{\mbox{\tiny ($m$-reg)}}$ be a regular $m$-gon with unit apothem (i.e.\ the shortest distance from its center to a side equals one). 
It is a pleasant exercise to show that $\rho_{m} = \myA(\mathcal{U}^{\mbox{\tiny ($m$-reg)}}) = \myS(\mathcal{U}^{\mbox{\tiny ($m$-reg)}})$. 
These facts together with Proposition 5.3 establish the following.

\noindent 
{\bf Corollary 5.4}\ \ {\sl Keep the notation of the preceding paragraph. 
If $\mathcal{M}$ is any unit $m$-gon, then} $\rho_{m} \leq \Pi_{\mathcal{M}}$, {\sl with equality if and only if} $\mathcal{M} \cong \mathcal{U}^{\mbox{\tiny ($m$-reg)}}$.\hfill\QED

{\bf \S 6\ \ Unital aspects of the unit shapes from our menagerie.}
There are many ways to find a unit distance within a given unit shape. 
For example, the unit circle is the only circle such that a central angle of $\pi/3$ subtends a chord of length $1$. 
But of course this hardly seems like the most natural way to distinguish the unit circle as a unital object. 

A more systematic idea is to locate a prominent point interior to 
a given shape and place the center of a unit circle at that point; then points of intersection between this unit circle and our given shape might be naturally distinguished. 
In fact, this approach is fruitful for unit right triangles, unit triangles in general, unit rhombi, and unit ellipses, as suggested by our figures below. 
But, this approach does not seem to be definitive. 
Indeed, part of what intrigues us about unit shapes is that there seems to be some art in discerning their interesting unital aspects. 

In the diagrams of Figure 6.1, we depict without much further explanation some unital aspects of the unit shapes from our \S 4 menagerie. 
For the unit right triangle and the generic unit triangle, we challenge the reader to discern the name of the famous central point depicted in each case (cf.\ \cite{Kimberling}). 
For the unit ellipse, the quantity $\myeulerR$ is given by $\displaystyle \myeulerR := \frac{1}{\pi}\int_{0}^{\pi}\sqrt{1+(\myr^{2}-1)\cos(t)^{2}\, }dt$. 

\begin{figure}[ht]
\begin{center}
\hspace*{0in}
\includegraphics[scale=1.05]{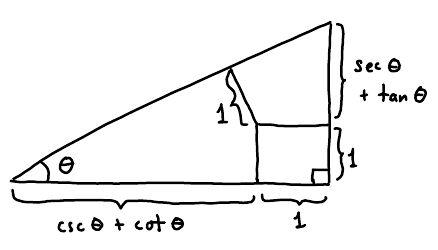} 
\hspace*{0in}
\includegraphics[scale=1.05]{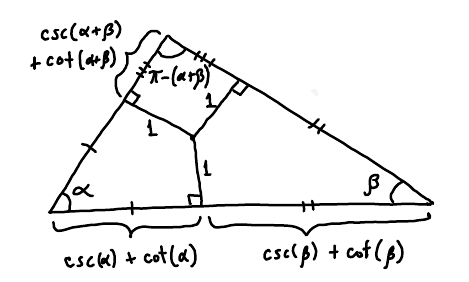} 
\hspace*{0in}
\includegraphics[scale=1.05]{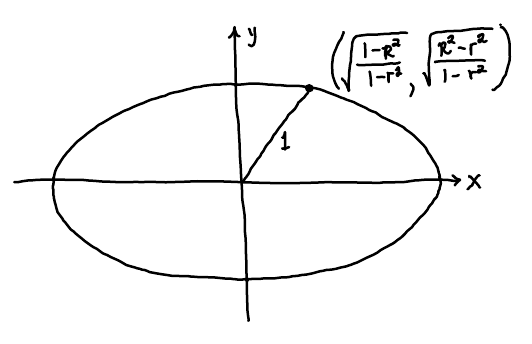} 
\end{center}
\begin{center}
\hspace*{-0.2in}
\includegraphics[scale=1.05]{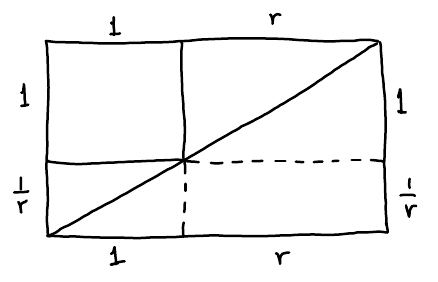} 
\hspace*{0in}
\includegraphics[scale=1.05]{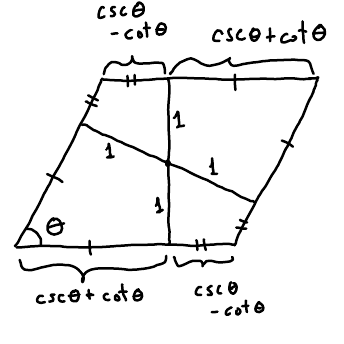} 
\hspace*{0in}
\includegraphics[scale=1.05]{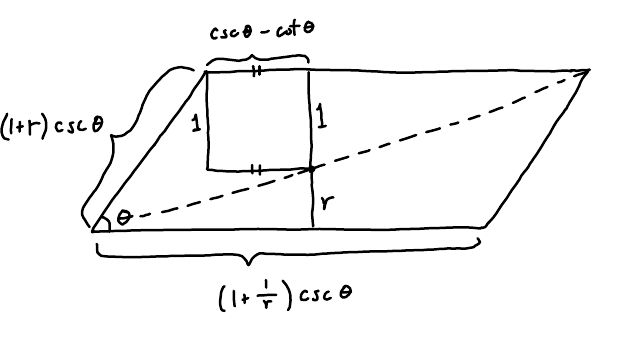} 

Figure 6.1: Unital aspects of some unit shapes.
\end{center}
\end{figure}

{\bf \S 7\ \ A comment on $\pi$ versus $\tau$.} 
In the {\em Scientific American} article \cite{Bartholomew}, Randyn Bartholomew writes about ``a vocal and growing minority of mathematicians who rally around the radical proposition that pi is wrong'' and that its alternative tau, where $\tau = 2\pi$, is the better constant to associate with circles and all $\pi$-related phenomena. 
Bartholomew's summary of the case tauists make for $\tau$ against $\pi$ is threefold: $\tau$ is the ratio of the circle's circumference to its radius, and a circle's radius is more important mathematically than its diameter; in its many and diverse appearances across the mathematical landscape, $\pi$ ``is preceded by a 2 more often than not'' and this seems to indicate that $\tau$ is the more intrinsic constant; and in pedagogy, fundamental concepts of trigonometry are ``distorted by this confusing factor of two.'' 

These tauistic assertions are not fully convincing to us, in part for reasons pertaining to the perspective of the present paper.  
We believe that the reflexive convention of regarding ``two times semiperimeter'' as ``perimeter'' in circumstances such as equation (2) from \S 1 is perhaps a perfunctory reading that devalues the importance of the factors ``two'' and ``semiperimeter'' as independent quantities. 
As we see it, semiperimeter seems to be a naturally occurring quantity in many expressions, such as Heron's formula for the area of a triangle as well as in our own rendering of the Isoperimetric Inequality in Theorem 5.1. 
Certainly as a fundamental measure for unit shapes, semiperimeter seems to be the more apt linear companion to area than perimeter. 
Also, consider that for a calculus-friendly indexing $\{\mathcal{C}_{\lambda}\}_{\lambda \in \mathbb{R}^{+}}$ of a family of shapes similar to a unit shape $\mathcal{U} := \mathcal{C}_{1}$ with fundamental measure $\Pi_{\mathcal{U}}$, we have $\myA(\lambda) = \lambda^{2}\Pi_{\mathcal{U}}$, $\myS(\lambda) = \lambda\Pi_{\mathcal{U}}$, and $\myA'(\lambda) = \frac{d}{d\lambda}(\lambda^{2}\Pi_{\mathcal{U}}) = 2\lambda\Pi_{\mathcal{U}} = 2\myS(\lambda)$. 
The factor of ``$2$'' in the expression $\myA'(\lambda) = 2\myS(\lambda)$ seems to us to be best understood as a power rule coefficient rather than a semiperimeter-doubling coefficient. 
Or, consider a generalization of the differential calculations of the first paragraph of \S 1. 
We have 
\[\Delta\myA = (\lambda+d\lambda)^{2}\Pi_{\mathcal{U}}-\lambda^{2}\Pi_{\mathcal{U}} = 2\left[\frac{\lambda+(\lambda+d\lambda)}{2}\right]d\lambda\,  \Pi_{\mathcal{U}}.\] 
In the latter expression, the factor of ``$2$'' plays a distinct role in allowing us to think of the area of a strip along the boundary precisely as the width of the strip multiplied by twice the average Tong inradius multiplied by the fundamental measure $\Pi_{\mathcal{U}}$. 

For circles, we agree that the radius of a circle is its crucial parameter in part because radius is a calculus-friendly indexing parameter for the family circles. 
Within such a context, perhaps we should be thinking of $\pi$ not only as the ratio of a circle's area to the square of its radius and as the area of a unit circle but also as the ratio of a circle's semiperimeter to its radius and as the semiperimeter of a unit circle. 
Perhaps the ``confusing factor of two'' should not be viewed so much as a party to the constant $\pi$ but rather as a distinguished modifier withinin the expression $\pi \cdot 2r$. 

{\bf \S 8\ \ Some other possible problems.} 

\noindent 
{\bf Problem 8.1 (open)}\ \ Add to the menagerie of unit shapes in \S 4 by investigating unit isosceles trapezoids, unit Norman windows, unit Star Trek communicator badges (formed via parabolic arcs), etc.  

\hspace*{-0.1in}
\includegraphics[scale=1.2]{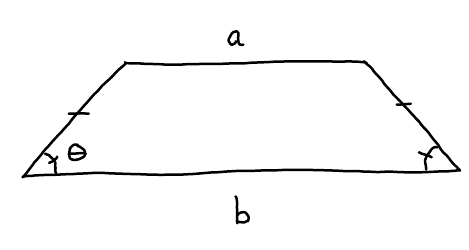} 
\hspace*{0.2in}
\includegraphics[scale=1.2]{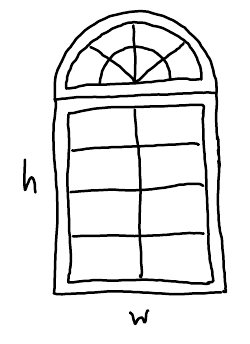} 
\hspace*{0.2in}
\includegraphics[scale=1.2]{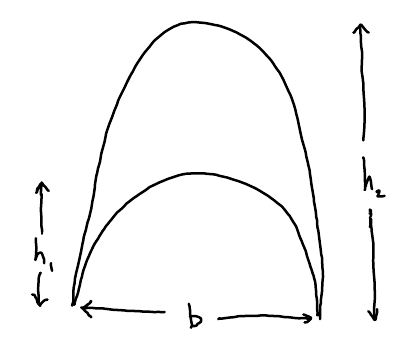} 

\noindent 
{\bf Problems 8.2}\ \  Solve the following problems without evaluating the integrals in question and without invoking any trigonometry.  (a) Show that the semiperimeter of the unit circle centered at the origin is the improper integral $\displaystyle \int_{-1}^{1}\frac{1}{\sqrt{1-x^{2}\, }}dx$. {\sc note:} We found it helpful to use the parameterization $t \mapsto \left(\frac{2t}{1+t^{2}},\frac{1-t^{2}}{1+t^{2}}\right)$, with $-1 \leq t \leq 1$, for the upper half of the unit circle.   
(b) Prove that the unit circle is a unit shape by showing that its area $\displaystyle \int_{-1}^{1}2\sqrt{1-x^{2}\, }dx$ equals its semiperimeter $\displaystyle \int_{-1}^{1}\frac{1}{\sqrt{1-x^{2}\, }}dx$. 

\noindent 
{\bf Problem 8.3 (open)}\ \ We have observed that Theorems 5.1 and 5.2 are equivalent. 
In view of this equivalence, can one produce a new proof of the Isoperimetric Inequality by proving Theorem 5.2 directly and then deducing Theorem 5.1 as a corollary? 

\noindent 
{\bf Problem 8.4 (open)}\ \ Generalize the unit shape notion to three (or more) dimensions, and produce some interesting examples. 
For example, for solids in $\mathbb{R}^{3}$, we believe that a unit shape should have the property that its volume is one-third its surface area. 
With this understanding, the unit platonic solids are those platonic solids whose inscribed sphere has a unit length radius. 
Then the fundamental measures of the platonic solids are: 
\begin{center}
\begin{tabular}{|c||c|c|c|c|c|}
\hline 
Platonic solid & Tetrahedron & Cube & Octahedron & Dodecahedron & Icosahedron\\
\hline
Fundamental measure & $8\sqrt{3}$ & $8$ & $4\sqrt{3}$ & \rule[-5.5mm]{0mm}{13mm}$\displaystyle \frac{20\xi}{\varphi^{3}}$ & \rule[-5.5mm]{0mm}{13mm}$\displaystyle \frac{20\sqrt{3}}{\varphi^{4}}$\\
\hline 
\end{tabular}

$\left(\mbox{In this table, $\displaystyle \varphi := 2\cos \frac{\pi}{5}$ is the golden ratio, and $\displaystyle \xi := 2\sin \frac{\pi}{5}$.}\right)$ 
\end{center}

\renewcommand{\refname}{\normalsize \bf References}
\renewcommand{\baselinestretch}{1}

\end{document}